\documentclass{article}

\input{diagrams}
\diagramstyle[repositionpullbacks,midshaft,nohug,scriptlabels]

\newcommand{\mcm}[3]{\newcommand{#1}[#2]{{\ensuremath{#3}}}}

\usepackage{latexsym}
\usepackage{amssymb}

\mcm{\emptybk}{0}{\:\:}
\mcm{\blank}{0}{(\emptybk)}
\mcm{\dashbk}{0}{-}
\mcm{\hyph}{0}{\mbox{-}}

\mcm{\diagspace}{0}{\mbox{\hspace{2em}}}

\newcommand{\bref}[1]{(\ref{#1})}

\mcm{\mb}{1}{\mathbf{#1}}
\mcm{\mc}{1}{\mathcal{#1}}
\mcm{\mi}{1}{\mathit{#1}}
\mcm{\mr}{1}{\mathrm{#1}}
\mcm{\cat}{1}{\mc{#1}}
\mcm{\fcat}{1}{\mb{#1}}
\mcm{\ovln}{1}{\overline{#1}}
\mcm{\twid}{1}{\widetilde{#1}}

\newcommand{\slsh}{/\linebreak[0]}

\newcommand{\dt}{.\linebreak[0]}

\mcm{\sub}{0}{\ \subseteq\ }
\mcm{\such}{0}{\:|\:}
\mcm{\without}{0}{\setminus}

\mcm{\ladj}{0}{\,\dashv\,}
\mcm{\of}{0}{\raisebox{0.08ex}{\ensuremath{\scriptstyle\circ}}}
\mcm{\sof}{0}{\raisebox{0.08ex}{\ensuremath{\scriptscriptstyle\circ}}}

\mcm{\bdry}{0}{\partial}
\mcm{\blob}{0}{\raisebox{.3ex}{\ensuremath{\scriptscriptstyle{\bullet}}}}

\mcm{\implies}{0}{\,\Rightarrow\,}

\mcm{\Hom}{0}{\mr{Hom}}
\mcm{\ob}{0}{\mr{ob}\,}		
\mcm{\op}{0}{\mr{op}}
\mcm{\comp}{0}{\mr{comp}}
\mcm{\id}{0}{\mr{id}}
\mcm{\ids}{0}{\mr{ids}}

\mcm{\Ab}{0}{\fcat{Ab}}
\mcm{\Alg}{0}{\fcat{Alg}}
\mcm{\Bicat}{0}{\fcat{Bicat}}
\mcm{\Bim}{1}{\fcat{Bim}(#1)}
\mcm{\Cat}{0}{\fcat{Cat}}
\mcm{\fc}{0}{\fcat{fc}}
\mcm{\Gph}{0}{\fcat{Gph}}
\mcm{\Graph}{0}{\fcat{Graph}}
\mcm{\Multicat}{0}{\fcat{Multicat}}
\mcm{\One}{0}{\fcat{1}}
\mcm{\Set}{0}{\fcat{Set}}
\mcm{\Span}{0}{\fcat{Span}}
\mcm{\Struc}{0}{\fcat{Struc}}
\mcm{\Sym}{0}{\fcat{Sym}}
\mcm{\Top}{0}{\fcat{Top}}
\mcm{\UBicat}{0}{\fcat{UBicat}}

\mcm{\integers}{0}{\mathbb{Z}}

\mcm{\range}{2}{#1,\,\ldots\,,#2}
\mcm{\tuplebts}{1}{(#1)}
\mcm{\bftuple}{2}{\tuplebts{\range{#1}{#2}}}
\mcm{\tuple}{3}{\tuplebts{\range{#1,#2}{#3}}}

\mcm{\eend}{2}{#1[#2]}
\mcm{\ehom}{3}{#1[#2,#3]}
\mcm{\ftrcat}{2}{[#1,#2]}

\mcm{\goesto}{0}{\,\longmapsto\,}
\mcm{\goiso}{0}{\goby{\diso}}
\mcm{\monic}{0}{\rMonic}
\mcm{\og}{0}{\lTo}
\mcm{\ogby}{1}{\lTo^{#1}}

\mcm{\oppair}{2}{\pile{\rTo^{\scriptstyle #1}\\ \lTo_{\scriptstyle #2}}}
\mcm{\parpair}{2}{\pile{\rTo^{\scriptstyle #1}\\ \rTo_{\scriptstyle #2}}}
\mcm{\parpairu}{0}{\pile{\rTo\\ \rTo}}

\mcm{\vslob}{3}
	{\left.
	\begin{diagram}[height=1.5em]
	#1		\\
	\dTo>{\,#2}	\\
	#3		\\
	\end{diagram}
	\right.}

\newarrow{Equals}=====
\newarrow{Get}....>
\newarrow{Goesto}|---{->}
\newarrow{Incl}C--->
\newarrow{Mod}--+->
\newarrow{Monic}{vee}---{vee}
\newarrow{NT}===={=>}

\mcm{\node}{0}{\bullet}
\mcm{\enode}{0}{\circ}
\mcm{\nl}{1}{\stackrel{\textstyle #1}{\node}}

\mcm{\diso}{0}{\sim}
\mcm{\vdiso}{0}{\wr}

\newcommand{\go}{\rTo\linebreak[0]}
\newcommand{\goby}[1]{\rTo^{#1}\linebreak[0]}
\newcommand{\nat}{\mathbb{N}}	
\newcommand{\iso}{\cong}
\newcommand{\demph}[1]{\textbf{\textup{#1}}}
\newcommand{\done}{\hfill\ensuremath{\Box}}
\newenvironment{prooflike}[1]{\begin{trivlist}\item\textbf{#1}\ }
{\end{trivlist}}
\newenvironment{proof}{\begin{prooflike}{Proof}}{\end{prooflike}}

\def\today{\number\day\space \ifcase\month\or
  January\or February\or March\or April\or May\or June\or
  July\or August\or September\or October\or November\or December\fi
  \space\number\year}

\newcommand{\rev}[1]{#1^\mr{rev}}
\newcommand{\revs}[2]{#1^\mr{rev}_{#2}}
\newcommand{\ofrev}{\of_\mr{rev}}
\newcommand{\ho}{{[0,1)}}
\newcommand{\osprod}{\boxplus}

\newtheorem{lemma}{Lemma}

\title{Are operads algebraic theories?}
\author{Tom Leinster}
\date{\normalsize  
University of Glasgow\\
T.Leinster\mbox{@}maths.gla.ac.uk\\
www.maths.gla.ac.uk/$\sim$tl}

\begin{document}

\maketitle

\begin{abstract}
I exhibit a pair of non-symmetric operads that, although not themselves
isomorphic, induce isomorphic monads.  The existence of such a pair implies
that if `algebraic theory' is understood as meaning `monad', operads cannot
be regarded as algebraic theories of a special kind.
\end{abstract}

\section*{Introduction}

Operads tend to be thought of as algebraic theories of some kind,
with the $n$th piece $P(n)$ of an operad $P$ thought of as the
collection of $n$-ary operations.  This point of view seems to be
validated by the fact that any operad has a category of algebras.
Nevertheless, it is not clear in principle that the passage from
an operad to its category of algebras does not involve a loss of
information.  The purpose of this note is to show by example that
such a loss can indeed occur, in the setting of operads without
symmetric group action.

The passage from operads to algebraic theories can be expressed
more precisely as follows.  I use `operad' to mean `non-symmetric
operad of sets'.  Any operad induces a monad on $\Set$, the
algebras for which are exactly the algebras for the operad.  Any
map of operads induces a map between the resulting monads, where
by definition a map of monads is a natural transformation
preserving multiplication and units in an obvious sense made
precise below.  This defines a functor
\[
(\textrm{operads}) \go (\textrm{monads on }\Set).
\]
But this functor does not reflect isomorphism: in other words,
there exist non-isomorphic operads $P$ and $P'$ whose associated
monads \emph{are} isomorphic.  This implies, of course, that the
categories of algebras for $P$ and $P'$ are isomorphic, so $P$
and $P'$ are `Morita equivalent' in a strong sense.  It also
implies that an operad should not be regarded as merely a monad
with certain properties: the canonical map from isomorphism
classes of operads to isomorphism classes of monads is not
injective.

Such a pair of operads is constructed as follows.  Any operad $P$ gives
rise to a new operad $\rev{P}$, whose induced monad is isomorphic to that
of $P$ (Section~\ref{sec:rev}).  It is then just a matter of finding an
operad $P$ such that $P \not\iso \rev{P}$.
This is done in Section~\ref{sec:counter}; further
comments follow in Section~\ref{sec:further}.

\paragraph*{Acknowledgements}
I thank Michael Batanin, Steve Lack, Peter May and Nathalie Wahl for their
remarks.

\section{The reverse of an operad}
\label{sec:rev}

For each operad $P$, I define its `reverse' $\rev{P}$ and show
that the monads induced by $P$ and $\rev{P}$ are isomorphic.

Let $P$ be an operad.  Its \demph{reverse} $\rev{P}$ is defined as follows:
$\rev{P}(n) = P(n)$ for all $n\in\nat$, the identity of $\rev{P}$ is the
same as that of $P$, and the composition $\ofrev$ is given by
\[
\theta \ofrev (\theta_1, \ldots, \theta_n)
=
\theta \of (\theta_n, \ldots, \theta_1)
\]
($n, k_i \in \nat$, $\theta \in P(n)$, $\theta_i \in P(k_i)$).  This
does define an operad $\rev{P}$: all that needs checking is associativity,
which is straightforward.

Let $(S, \mu, \eta)$ be the monad on $\Set$ induced by $P$.  Then for any
set $X$,
\[
SX =
\sum_{n\in\nat} P(n) \times X^n,
\]
the unit map
\[
\eta_X: X \go SX
\]
picks out the identity element of $P(1)$, and the multiplication map
\[
\mu_X: S^2 X \go SX
\]
is given by
\begin{eqnarray*}
&&\left(\theta, 
(\theta_1, x_1^1, \ldots, x_1^{k_1}),
\ldots,
(\theta_n, x_n^1, \ldots, x_n^{k_n})\right)\\
&\goesto &
(\theta \of (\theta_1, \ldots, \theta_n),
x_1^1, \ldots, x_1^{k_1}, \ldots, x_n^1, \ldots, x_n^{k_n})
\end{eqnarray*}
($n, k_i \in \nat$, $\theta \in P(n)$, $\theta_i \in P(k_i)$, $x_i^j \in X$).
The monad $(\rev{S}, \rev{\mu}, \rev{\eta})$ induced by $\rev{P}$ is the
same except that the multiplication formula becomes
\begin{eqnarray*}
&&\left(\theta, 
(\theta_1, x_1^1, \ldots, x_1^{k_1}),
\ldots,
(\theta_n, x_n^1, \ldots, x_n^{k_n})\right)\\
&\goesto &
(\theta \of (\theta_n, \ldots, \theta_1),
x_1^1, \ldots, x_1^{k_1}, \ldots, x_n^1, \ldots, x_n^{k_n}).
\end{eqnarray*}

There is a natural isomorphism $\iota: S \goiso \rev{S}$ whose component at
a set $X$ is
\[
\begin{array}{rrcl}
\iota_X:	&SX	&\goiso	&\rev{S} X	\\
		&
(\theta, x_1, \ldots, x_n)	&
\goesto	&
(\theta, x_n, \ldots, x_1)
\end{array}
\]
($n\in\nat$, $\theta\in P(n)$, $x_i \in X$).  Using the above descriptions 
of the monad structures, it is straightforward to check that $\iota$ is an
isomorphism of monads, in other words, that for each $X$ the diagrams
\[
\begin{diagram}[size=2em]
X		&\rEquals	&X	\\
\dTo<{\eta_X}	&		&\dTo>{\revs{\eta}{X}}	\\
SX		&\rTo_{\iota_X}	&\rev{S} X	\\
\end{diagram}
\diagspace
\begin{diagram}[size=2em]
S^2 X			&\rTo^{S \iota_{X}}	&S \rev{S} X	&
\rTo^{\iota_{\rev{S} X}}	&(\rev{S})^2 X		\\
\dTo<{\mu_X}		&			&		&
				&\dTo>{\revs{\mu}{X}}	\\
S X			&			&\rTo_{\iota_X}	&
				&\rev{S} X		\\
\end{diagram}
\]
commute.

So, as promised, any operad $P$ gives rise to a new operad $\rev{P}$
inducing the same monad as $P$.  

There are at least two abstract perspectives on this
construction.  First, write $T$ for the free monoid monad on
$\Set$.  Then an operad amounts to a cartesian monad $S = (S,
\mu, \eta)$ on $\Set$ together with a cartesian natural
transformation $\pi: S \go T$ respecting the monad structures,
and the monad induced by the operad is simply $S$.  (For
explanation and proof, see for instance Cor~6.2.4
of~\cite{HOHC}.)  Now, there is an involution $\rho$ of the monad
$T$ given by reversing the order of finite lists, which implies
that any operad $P$ described by a pair $(S, \pi)$ gives rise to
a new operad described by the pair $(S, \rho\of\pi)$; this operad
is $\rev{P}$.  From this point of view, the monad induced by
$\rev{P}$ is not just isomorphic but equal to that induced by
$P$.

Second, given any cartesian monad $T$ on a category $\cat{E}$
with pullbacks, there is a category of so-called $T$-operads.
(See Chapter~4 of~\cite{HOHC}.)  Any $T$-operad $P$ induces a
monad $T_P$ on $\cat{E}$, and algebras for the operad are by
definition algebras for this monad.  When $T$ is the free monoid
monad on $\Set$, these are the standard notions of non-symmetric
operad, induced monad, and algebra.  Inevitably, if we have an
isomorphism $(\cat{E}, T) \goiso (\cat{E}', T')$ between two
different cartesian monads on two different categories then there
is an induced isomorphism between the categories of $T$-operads
and $T'$-operads, and if $P'$ is the $T'$-operad corresponding to
a $T$-operad $P$ then the monad $T'_{P'}$ on $\cat{E'}$
is obtained by transporting the monad $T_P$ on $\cat{E}$
across the isomorphism.  In particular, this holds for the isomorphism
\[
(\id, \rho): (\Set, T) \goiso (\Set, T)
\]
where $T$ is the free monoid monad and $\rho$ is as above; the resulting
automorphism of the category of non-symmetric operads sends $P$ to
$\rev{P}$, and by the preceding comments the respective induced monads are
isomorphic.

Observe also that reversal works for (non-symmetric) operads in any
symmetric monoidal category $\cat{V}$.  The definition of $\rev{P}$ is an
absolutely straightforward generalization of the case $\cat{V} = \Set$,
using the symmetry of $\cat{V}$.  If $\cat{V}$ has countable coproducts and
tensor distributes over them then any operad $P$ in $\cat{V}$ induces a
monad on $\cat{V}$, algebras for which are algebras for $P$; and just as
above, the monads induced by $P$ and $\rev{P}$ are isomorphic.

\section{The counterexample}
\label{sec:counter}

To find a pair of non-isomorphic operads whose induced monads are
isomorphic, it suffices to find an operad not isomorphic to its reverse.  

This task is not completely straightforward, since many commonly
encountered operads admit a symmetric structure and any such operad
\emph{is} isomorphic to its reverse.  Indeed, let $\sigma_n \in S_n$ denote
the permutation reversing the order of $n$ letters: then for any symmetric
operad $P$, there is an isomorphism $P \goiso \rev{P}$ sending $\theta \in
P(n)$ to $\theta \cdot \sigma_n \in \rev{P}(n)$.  Further, several
well-known operads that do not admit a symmetric structure are,
nevertheless, isomorphic to their reverse: this applies, for instance, to
Stasheff's operad of associahedra (\cite{StaHAHI}, \cite{MSS}).

Here is an operad $P$ not isomorphic to its reverse.  Let $P(n)$ be the set
of all $n$-tuples $(f_1, \ldots, f_n)$ of order-preserving continuous maps
$f_i: \ho \go \ho$ of the half-open real interval such that if $i<j$ and
$t_i, t_j \in \ho$ then $f_i(t_i) < f_j(t_j)$.  The identity of $P$ is the
identity map $\id_\ho \in P(1)$, and composition
\[
P(n) \times P(k_1) \times \cdots \times P(k_n)
\go
P(k_1 + \cdots + k_n)
\]
is 
\begin{eqnarray*}
&&
((f_1, \ldots, f_n),
(f_1^1, \ldots, f_1^{k_1}), \ldots, (f_n^1, \ldots, f_n^{k_n}))
\\
&\goesto &
(f_1 f_1^1, \ldots, f_1 f_1^{k_1}, \ldots, 
f_n f_n^1, \ldots, f_n f_n^{k_n}).
\end{eqnarray*}
Seen another way, $P$ is an endomorphism operad.  For consider the
(non-symmetric) monoidal category of ordered topological spaces, where the
product $X \osprod Y$ is defined by taking the disjoint union of $X$ and
$Y$ and adjoining the relation $x < y$ for each $x \in X$ and $y \in Y$.
(Compare addition of ordinals.)  Then $P(n)$ is the set of maps
$\ho^{\osprod n} \go \ho$, with the usual endomorphism operad structure.

To prove that $P$ is not isomorphic to $\rev{P}$, I introduce some
temporary terminology.  Let $Q$ be an operad.  An element $\gamma\in Q(1)$
is \demph{constant} if
\[
\textrm{for all $n\in\nat$ and all $\phi, \phi' \in Q(n)$,\ }
\gamma \of (\phi) = \gamma \of (\phi').
\]
An element $\phi \in Q(n)$ is \demph{surjective} if 
\[
\textrm{for all $\theta, \theta' \in Q(1)$,\ }
\theta \of (\phi) = \theta' \of (\phi)
\,\implies\, 
\theta = \theta'.
\]

The following lemma shows that these terms have the expected meanings when
$Q$ is $P$ or $\rev{P}$.  For convenience, I write an element $(g) \in
P(1)$ as simply $g$.
\begin{lemma} \ 
\begin{enumerate}
\item \label{constant}
  $g \in P(1)$ is constant in the sense above if and only if the map
  $g: \ho \go \ho$ is constant in the usual sense.
\item \label{surjective} 
  $(f_1, \ldots, f_n) \in P(n)$ is surjective in the sense above if and
  only if the union of the images of $f_1, \ldots, f_n$ is $\ho$.
\end{enumerate}
Moreover, both statements remain true when $P$ is replaced by $\rev{P}$.
\end{lemma}

\begin{proof}
For~\bref{constant}, `if' is clear.  Now suppose $g$ is not constant in the
usual sense, so that there exist $t, t' \in \ho$ with $g(t) \neq g(t')$.
If we take $f, f': \ho \go \ho$ to be the constant functions with
respective values $t$ and $t'$ then $f, f' \in P(1)$ with $g \of (f) \neq g
\of (f')$, so $g$ is not constant in the sense above.

For~\bref{surjective}, `if' is also clear.  Conversely, if the union of the
images of $f_1, \ldots, f_n$ is not $\ho$ then by continuity, one of the
following holds:
\begin{itemize}
\item $n = 0$
\item $n\geq 1$ and $f_1(0) > 0$
\item $\sup f_{i-1} < f_i(0)$ for some $i \in \{2, \ldots, n\}$
\item $n\geq 1$ and $\sup f_n < 1$. 
\end{itemize}
In all cases, there is some nonempty open interval $(a, b) \sub \ho$
that does not meet the union of the images of $f_1, \ldots, f_n$.  We can
construct a continuous order-preserving map $h: \ho \go \ho$ that is not
the identity but satisfies $h(t) = t$ for all $t \not\in (a, b)$,
and this gives distinct elements $h, \id$ of $P(1)$ satisfying $h \of (f_1,
\ldots, f_n) = \id \of (f_1, \ldots, f_n)$.  So $(f_1, \ldots, f_n)$ is not
surjective in the sense above.

`Moreover' follows immediately from the definition of $\rev{P}$. \done
\end{proof}

We can now show that the following isomorphism-invariant property of an
operad $Q$ holds when $Q = P$ but fails when $Q = \rev{P}$:
\begin{quote}
  there exist $\phi \in Q(2)$ and constant $\gamma \in Q(1)$ such that
  $\phi \of (\gamma, \id) \in Q(2)$ is surjective.
\end{quote}
It will follow that $P \not\iso \rev{P}$.

To see that the property holds for $Q = P$, let $g, f_1: \ho \go \ho$ both
be the map with constant value $0$, and let $f_2: \ho \go \ho$ be the
identity.  Then $g \in P(1)$ is constant, $\phi = (f_1, f_2)$ is an element
of $P(2)$, and if $e = \phi \of (g, \id)$ then
\[
e_2 = f_2 \of \id = \id
\]
so $e$ is surjective.

To see that the property fails for $Q = \rev{P}$, we have to see that given
$\phi = (f_1, f_2) \in P(2)$ and constant $g \in P(1)$, the composite $e =
\phi \of (\id, g)$ in $P$ cannot be surjective.  Indeed, let $b$ be the
constant value of $g$: then 
\[
e_2 = f_2 \of g = \textrm{(constant map with value $f_2(b)$)},
\]
so by order-preservation
\[
\textrm{image}(e_1) \cup \textrm{image}(e_2)
\ \subseteq\ 
[0, f_2(b)]
\ \subsetneqq\ 
\ho,
\]
as required.

\section{Further comments}
\label{sec:further}

The properties of the functor 
\[
G: (\textrm{operads}) \go (\textrm{monads on }\Set)
\]
can be analyzed more precisely.  The monads in the essential
image of $G$ (that is, the monads isomorphic to $G(P)$ for some
operad $P$) are the strongly regular finitary monads.  By
definition, this is the class of monads whose corresponding
algebraic theory can be presented by finitary operations and
equations in which the same variables appear on each side of the
equals sign, in the same order and without repetition.  (For
instance, the theory of monoids is allowed, but the theories of
commutative monoids and groups are not.)  Another description is
that they are the cartesian monads $S$ such that there exists a
cartesian natural transformation, respecting the monad
structures, from $S$ to the free monoid monad.  The original
source on strong regularity is~\cite{CJ}; proofs of the results
just mentioned can be found in~C.1 and~6.2 of~\cite{HOHC}.

The functor $G$ does not reflect isomorphism, as has been shown.
It does reflect isomorphisms (plural): that is, if $f: P \go P'$
is a map of operads and $G(f)$ is an isomorphism then so too is
$f$.  This is easily shown, as is the fact that $G$ is faithful.
But since a full and faithful functor reflects isomorphism, $G$
cannot be full.  To prove this more directly, let $P$ be the functor
of the previous section and take the isomorphism $\iota: S \goiso
\rev{S}$ of Section~\ref{sec:rev}, where $S$ and $\rev{S}$ are
the monads induced by $P$ and $\rev{P}$ respectively.  Then since
$G$ reflects isomorphisms, there is no map $f: P \go \rev{P}$
satisfying $G(f) = \iota$: so again, $G$ is not full.

Here I have stuck to operads of sets; I know little about the
situation for other types of operad.  Trivially, taking discrete
spaces on the set-theoretic example above yields a pair of
non-symmetric topological operads that induce isomorphic monads
but are not themselves isomorphic.

The situation for symmetric operads is completely different:
symmetric operads of sets \emph{can} be identified as monads of a
special kind.  Precisely, the canonical functor
\[
(\textrm{symmetric operads}) \go (\textrm{monads on }\Set)
\]
defines an equivalence between the category of symmetric operads
and the category of analytic monads and weakly cartesian maps.
This is a result of 
Weber~\cite{Web}, 
using Joyal's
characterization of the endofunctors on $\Set$ induced
by species and of the natural transformations induced by maps
between them~\cite{Joy}.

\end{document}